\documentstyle{amsart}
\input{bull-art}
\theoremstyle{plain}
\newtheorem{thm}{Theorem}

\theoremstyle{remark}
\newtheorem{rem}{Remark}

\newcommand{\App}{\operatorname{App}}
\newcommand{\avg}{\operatorname{avg}}
\newcommand{\reals}{\Bbb{R}}
\newcommand{\comp}{\operatorname{comp}}
\newcommand{\cost}{{\mathrm{cos}}t}
\newcommand{\Int}{\operatorname{Int}}

\begin{document}
\def\currentvolume{28}
\def\currentissue{2}
\def\currentyear{1993}
\def\currentmonth{April}
\def\copyrightyear{1993}
\def\currentpages{308-314}
\ratitle
\title[Integration and Approximation of Multivariate 
Functions]{Integration And Approximation\\ 
    of Multivariate Functions:\\
    Average Case Complexity with\\ 
    Isotropic Wiener Measure}
\author{G. W. Wasilkowski}
\thanks{This research was supported in part 
 by the National Science Foundation under Grant 
CCR-91-14042}
\address{Department of Computer Science, University of 
Kentucky, Lexington,
Kentucky 40506}
\email{greg@@ms.uky.edu}
\date{December 19, 1991}
\maketitle
\begin{abstract}
We study the average case complexity of multivariate 
integration and $L_2$ function approximation for the class 
$F=C([0,1]^d)$ of continuous functions of $d$ variables. 
The class $F$ 
is endowed with the isotropic Wiener measure (Brownian 
motion in Levy's sense).
Furthermore, for both problems, only function values are 
used as data. 
\end{abstract}

\begin{center} {\sc 1. Introduction}\end{center}
We study the integration and function approximation 
problems for multivariate 
functions $f$. For the integration problem, we want to 
approximate 
the integral of $f$ to within a specified error 
$\varepsilon$; and for the 
function approximation problem, we want to recover $f$ 
with the $L_2$ error 
not exceeding $\varepsilon$. To solve both problems, we 
would like to use as 
small a number of function values as possible. 

Both problems have been extensively studied in the 
literature (see, 
[9, 16] for hundreds of references). However, they are 
mainly 
addressed in the worst-case setting. In the worst-case 
setting the cost and 
the error of an algorithm are defined by the worst 
performance with 
respect to the given class $F$ of functions $f$. Not 
surprisingly, 
for a number of classes $F$, the integration and function 
approximation 
problems are intractable (prohibitively expensive) or even 
unsolvable. 
For instance, if $F$ consists of continuous functions that 
are bounded by 1, 
no algorithm that uses a finite number of function 
values can approximate the integral of $f$, nor can it 
recover $f$ with the 
worst-case error less than 1. Hence, both problems are 
unsolvable for 
$\varepsilon<1$. Assuming that functions $f$ have bounded 
$r$th derivative in 
the sup-norm, the number of function values required for 
the worst-case 
error not to exceed $\varepsilon$ is of order 
$\varepsilon^{-d/r}$. 
Hence, for fixed $r$, it is exponential in $d$. 

Due to intractability in the worst-case setting, the 
average-case setting 
is of interest. In the average-case setting, the class $F$ 
is 
equipped with a probability measure $\mu$. The error and 
the cost of an 
algorithm are measured by the expectations with respect to 
$\mu$. Then, 
the average-case complexity (with respect to $\mu$) is 
defined as the minimal 
expected cost needed to compute an approximation with the 
expected error not greater than $\varepsilon$. 

The majority of the average-case results obtained so far 
(see, 
[3--6, 9--18, 21]) 
deal with scalar functions ($d=1$). These results indicate 
that for a 
``reasonable'' choice of measure $\mu$, the integration 
and function 
approximation problems are significantly easier on the 
average than in the 
worst-case setting. Thus, one could hope that the 
intractability (or even 
noncomputability) of multivariate problems in the 
worst-case setting can be 
removed by switching from the worst-case to the 
average-case setting. 

This hope has recently been supported by Wo\'zniakowski  
(see 
[23, 24]), 
who analyzes integration and function approximation for 
the class $F=C([0,1]^d)$ endowed with the Wiener sheet 
measure $\mu$. 
He proves that the average-case complexities of both 
problems are only 
weakly dependent on the number of variables. Indeed, the 
average-case complexity of computing an 
$\varepsilon$-approximation is 
$\Theta(\varepsilon^{-1}(\log \varepsilon^{-1})^{(d-1)/2})$ 
for the integration problem and 
$\Theta(\varepsilon^{-2}(\log \varepsilon^{-1})^{2(d-1)})$ 
for the function approximation problem. 

In this paper we study the average-case complexity of 
the integration and function approximation problems. 
However, instead of 
the Wiener sheet measure, we endow the class 
$F=C([0,1]^d)$ with 
the isotropic Wiener measure (or Brownian motion in Levy's 
sense). 
We prove that the average-case complexity equals 
$\Theta(\varepsilon^{-2/(1+1/d)})$
for the integration problem and $\Theta(\varepsilon^{-2d})$ 
for the function approximation problem. Unlike the Wiener 
sheet measure, 
the average-case complexity of the function approximation 
problem depends 
strongly on $d$. In particular, for large $d$ this problem 
is intractable 
since its complexity 
$\Theta(\varepsilon^{-2d})$ is exponential in $d$ and is 
huge even for a modest 
error demand $\varepsilon$. For large $d$ the average-case 
complexity 
of the integration problem is essentially proportional to 
$\varepsilon^{-2}$, 
which is the highest possible average-case complexity of 
the integration 
problem. Indeed, for {\em any} probability measure with 
finite 
expected value of $\|f\|^2_{L_2}$, the average-case 
complexity is bounded 
from above by $O(\varepsilon^{-2})$. Hence, this is again 
a negative result. 

Thus, the average-case complexities of integration and 
function approximation 
problems are very different depending on whether $\mu$ is 
the Wiener sheet 
or isotropic Wiener measure. It is interesting to note 
that both measures 
are identical when $d=1$. They are different for $d>1$; 
results of 
[23, 24] and our results indicate how drastically 
different they are. 

The paper is organized as follows. Section 2 provides 
basic definitions.  The
main results are presented in \S 3. In addition to results 
already  mentioned,
\S 3 discusses optimality of Haber's [2]  modified Monte 
Carlo quadrature and
of a piecewise constant function  approximation. It also 
contains a result
relating the average-case complexities of the integration 
and function
approximation  problems for general probability measures. 
In this paper we omit
all  proofs because of their substantial length.

\vskip 1 pc
\begin{center}{\sc 2. Basic definitions}
\end{center}
In this paper  we consider the following integration and 
function 
approximation problems 
for multivariate functions. Let $F=C(D)$ be the space of 
functions 
$f:D\to \reals$ where $D$ is a bounded subset of 
$\reals^d$. For simplicity,
we take $D=[0,1]^d$ 
as a unit cube. For every $f\in F$ we wish to approximate 
$S(f)$, 
where $S:F\to G$ with 
\[
  S(f)=\mbox{Int}(f)=\int_Df(x)\,dx \quad \mbox{and}\quad 
G=\reals 
    \qquad\mbox{for the integration problem,}\]
\[
  S(f)=\mbox{App}(f)=f \quad \mbox{and}\quad G=L_2(D) 
    \qquad\mbox{for the approximation problem}.\]

We assume that the functions $f$ are unknown; instead we 
can compute 
information $N(f)$ that consists of a finite number of 
values of $f$ taken at some points from $D$. For a precise 
definition 
of $N$ see [16]. Here we stress only that 
\[
  N(f)=\left[f(x_1),\dots,f(x_n)\right],
\]
where the points $x_i$ and the number $n$ of them (called 
the 
{\em cardinality} of $N$) can be selected adaptively 
and/or randomly. 
That is, for adaptive $N$, $x_i$'s depend on previously 
computed values 
$f(x_1),\dots,f(x_{i-1})$, and the cardinality $n=n(f)$ 
varies with $f$ 
based on computed values. For randomized $N$, the points 
$x_i$ and the 
cardinality $n(f)$ may also depend on an outcome of a 
random process $t$. 
(That is, $x_i$ is selected randomly with an arbitrary 
distribution 
that may depend on previously computed values of $f$; the 
distribution of 
$n(f)$ may also depend on observed values.) 
In such a case, we sometimes write $N(f)=N_t(f)$. 

An approximation $U(f)$ to $S(f)$ is computed based on 
$N(f)$. That is, 
\[
   U(f)=\phi(N(f)), \quad\mbox{where}\ \phi:N(F)\to G
\]
is an arbitrary mapping; $\phi$ is called an {\em 
algorithm} that uses 
$N$. The algorithm $\phi$ can also be random; in such a 
case, we sometimes 
write $\phi=\phi_t$. 

In the average-case setting, we assume that the space $F$ 
is endowed with 
a (Borel) probability measure $\mu$. Then the {\em average 
error} 
and the {\em average cost}\footnote{We measure the cost by 
the 
expected number of function values neglecting the 
combinatory cost of $N$ 
and of $\phi$. With the exception of Theorem \ref{thm-2}, 
this is without loss 
of generality since, as explained in a number of 
references (see, e.g., 
[16]), for Gaussian measures the same results hold for a 
more general 
definition of the average cost, provided 
that a single arithmetic operation is no more expensive 
than 
a function evaluation.}  of $\phi$ are defined 
respectively by
\begin{align*}
   e^{\avg}(\phi,N,S,\mu):=&\sqrt{\mbox{E}_{\mu}\mbox{E}_t 
     (\|S(f)-\phi_t(N_t(f))\|_G^2)},\\
   \cost^{\avg}(\phi,N,S,\mu):=&\mbox{E}_{\mu}%
\mbox{E}_t(n(f)).   
\end{align*}
(By $\mbox{E}_{\mu}$ and $\mbox{E}_t$ we denote the 
expectations 
w.r.t.\,$\mu$ and $t$, respectively.) 
Of course, for {\em deterministic} $N$ and $\phi$, 
\begin{align*}
   e^{\avg}(\phi,N,S,\mu)=&\sqrt{\int_F\|S(f)-\phi(N(f))%
\|_G^2\,\mu(df)},\\
   \cost^{\avg}(\phi,N,S,\mu)=&\int_Fn(f)\,\mu(df).
\end{align*}

The {\em average-case complexity} is the minimal average 
cost for solving the 
problem to within a preassigned error accuracy 
$\varepsilon$. 
That is, 
\[
  \comp^{\avg}(\varepsilon,S,\mu):=\inf\left\{\cost^{%
\avg}(\phi,N,S,\mu)
   :e^{\avg}(\phi,N,S,\mu) \le\varepsilon\right\}.
\]
(We stress that the infimum above is taken with respect to 
all randomized 
$\phi$ and $N$.) 

In this paper we analyze the average-case complexity of 
the integration 
and function approximation problems ($S=\mbox{Int}$ and 
$S=\mbox{App}$) 
with the assumption that the 
probability $\mu$ is the {\em isotropic} Wiener measure. 
This measure is 
also referred to as the Brownian motion in Levy's sense. 
For more detailed 
discussion and properties of $\mu$ (see [1, 7, 8]). 
Here we only recall that $\mu$ is a zero-mean Gaussian 
measure with the 
correlation function
\[
  K(x,y)=\frac{\|x\|+\|y\|-\|x-y\|}2\quad\forall 
x,y\in\reals^d,
  \qquad \|x\|^2=\sum_{i=1}^dx_i^2.
\]

\vskip 1 pc
\begin{center}{\sc 3. Main results}\end{center}

\begin{thm}\label{thm-1}
For the integration and function approximation problems, 
\begin{align}\label{help-1}
\comp^{\avg}(\varepsilon,\mbox{\em Int},\mu)
   =&\Theta\left(\varepsilon^{-2/(1+1/d)}\right),\\
\comp^{\avg}(\varepsilon,\mbox{\em App},\mu)=&
\Theta\left(\varepsilon^{-2d}\right).
\end{align}
\end{thm}

For $d=1$, $\mu$ equals the classical Wiener measure. 
Hence, for scalar 
functions  this theorem follows from known results (see 
[12, 13, 19]). 

We now exhibit algorithms and information that are almost 
optimal. 
Let $n=p^d$. Partition $D$ into $n$ equal-sized cubes 
$U_i$, 
$U_i=x_i+[-1/(2p),+1/(2p)]^d$, each centered at $x_i$. 
For the integration problem, consider the 
following randomized information and algorithm due to 
Haber (see [2]): 
\begin{equation}\label{opt-inf-int} 
  N_n^{\Int}(f)=[f(t_1),\dots,f(t_n)]\quad\mbox{and}\quad 
  \phi_n^{\Int}(N_n^{\Int}(f))=\frac1n\sum_{j=1}^nf(t_j), 
\end{equation}
where $t_i$'s are uniformly distributed in $U_i$'s. 
For the function approximation problem, consider 
\begin{equation}\label{opt-inf-app} 
  N^{\App}_n(f)=[f(x_1),\dots,f(x_n)] \quad\mbox{and}\quad
\phi_n^{\App}(N^{\App}_n(f))=\sum_{i=1}^ng_i(\cdot)f(x_i) 
,
\end{equation} 
with $g_i$ being the indicator function for the set $U_i$. 

\begin{thm}\label{thm-4}
For every $n$, the average errors of $\phi_n^{\Int}$ and 
$\phi_n^{\App}$ are 
respectively equal to 
\[
  \frac{\sqrt{\int_{D}\int_{D}\|x-y\|/2\,dxdy}}{n^{1/2+
1/(2d)}}
  \qquad\mbox{and}\qquad
  \frac{\sqrt{\int_{D}\|x\|/2\,dx}}{n^{1/(2d)}}.
\]

These algorithms are almost optimal. 
Indeed, for $n^{\Int}(\varepsilon)$ and 
$n^{\App}(\varepsilon)$ given by 
\begin{align*}
  n^{\Int}(\varepsilon)= &
  \left\lceil 
\left(\varepsilon^{-2}\int_D\int_D\|x-y\|/2\,dxdy
     \right)^{1/(d+1)}\right\rceil ^{d},\\
  n^{\App}(\varepsilon)=&
  \left\lceil 
\varepsilon^{-2}\int_{D}\|x\|/2\,dx\right\rceil ^{d}, 
\end{align*}
$\phi_{n^{\Int}(\varepsilon)}^{\Int}$ and 
$\phi_{n^{\Int}(\varepsilon)}^{\App}$ have 
the average errors less than or equal to $\varepsilon$ and 
their 
costs are proportional to \RM{(\ref{help-1})} and 
\RM{(2)}, respectively. 
\end{thm}

\begin{rem}
In the worst-case setting with $F=C[0,1]^d$, Haber's 
modified Monte 
Carlo algorithm $\phi_n^{\Int}$ and the classical Monte 
Carlo algorithm 
$n^{-1}\sum_{i=1}^n 
f(t_i)$ (with $t_i$'s uniformly distributed in $D$) 
have (modulo constants) the same errors that 
are proportional to $1/\sqrt{n}$. It can be proven that 
the average 
error of the classical Monte Carlo algorithm equals 
$\sqrt{n^{-1}\int_D\int_D\|x-y\|/2\,dxdy}$.
Thus, it is 
precisely $n^{1/(2d)}$ times larger than the average error 
of $\phi_n^{\Int}$. 
\end{rem}

\begin{rem}
Although the information $N_n^{\Int}$ is randomized, the 
cardinality 
$n$ is fixed. Thus, the mean value theorem implies the 
existence of 
{\em deterministic} $N^*_n(f)=[f(x^*_1),\dots,f(x^*_n)]$ 
such that 
$\phi^*_n(N^*_n(f))=n^{-1}\sum_{i=1}^nf(x^*_i)$ has the 
average error 
not exceeding the average error of $\phi_n^{\Int}$. (We do 
not know 
the location of the points $x^*_i$; we only know that 
$x^*_i\in U_i$, 
$1\le i\le n$.) 
This and the fact that the algorithm and information given 
in 
(\ref{opt-inf-app}) are deterministic imply that 
randomization does not 
help for both problems. The lack of power of randomization 
holds for more 
general problems. Indeed, 
randomization does not help for linear $S$ and Gaussian 
$\mu$ (see [20, 22]). 
\end{rem}

The final theorem relates the average-case complexities of 
the 
integration and function approximation problems for an 
arbitrary 
(Borel) probability measure. 

\begin{thm}\label{thm-2}
Let $\nu$ be an arbitrary probability measure on $F$. If 
\[
   \comp^{\avg}(\varepsilon,{\Int},\nu)=\Omega\left(%
\varepsilon^{-p}\right)
\]
for some $p$ \RM(obviously, $p\le2$ whenever 
$\|f\|^2_{L_2(D)}$ has a finite 
$(\nu$-\RM) expectation\/\RM), then 
\[
   \comp^{\avg}(\varepsilon,{\App},\nu)
    =\Omega\left(\varepsilon^{-p(1+p/(2-p))}\right).
\]
\end{thm}

\begin{rem}
This theorem can easily be extended in a number of ways. 
For instance, 
it holds when the function approximation problem is 
considered with the 
$L_2(D)$-norm replaced by $\|f-f^*\|^2=\int_Dw(x)
(f(x)-f^*(x))^2\,dx$, a weighted norm, for some weight 
$w\ge 0$ 
and the integration 
problem defined by ${\Int}(f)=\int_Df(x)\sqrt{w(x)}\,dx$. 

It also holds in the {\em worst-case setting with 
randomization}. 
In this setting, instead of the expectation 
$\mbox{E}_{\nu}$, we take 
the supremum w.r.t.\,$f\in F_0$ ($F_0$ is a given subset 
of $F$) in the 
definitions of the error and cost (the expectation w.r.t.\ 
random $t$ 
remains). (For more detailed definitions, see [16]). 
Hence, again, in the worst-case setting with randomization, 
integration is an easier problem than  is the 
function approximation problem. We stress that this need 
not be true if the 
worst-case {\em deterministic} ({\em without 
randomization}) 
setting is considered, since for a number of classes $F_0$ 
the integration 
and approximation problems have asymptotically the same 
worst-case 
deterministic complexities (see, [9, 16]). 
\end{rem}

\begin{center}{\sc Acknowledgment}\end{center}
The author wishes to thank David R.\,Adams, Klaus Ritter, 
and 
Henryk  Wo\'zniakowski for valuable suggestions. 


\end{document}